\documentclass[12pt,thmsa]{article}
\usepackage{amsfonts}
\usepackage{amssymb}
\newtheorem{teo}{Theorem}[section]

\newtheorem{prop}[teo]{Proposition}
\newtheorem{defin}[teo]{Def\mbox{}inition}
\newtheorem{obs2}[teo]{Remark}

\newtheorem{tea}{Theorem}[subsection]

\newtheorem{no2}[teo]{Note}

\newtheorem{no3}[tea]{Note}

\newcommand{\Gal}{{\rm Gal}}
\newcommand{\Frob}{{\rm Frob }}
\newcommand{\trace}{{\rm trace}}
\newcommand{\mod}{{\rm mod}}

\newcommand{\Q}{\mathbb{Q}}

\newcommand{\PGL}{{\rm PGL}}

\newcommand{\GL}{{\rm GL}}

\newcommand{\cond}{{\rm cond}}

\newcommand{\eq}{x^4 + y^4 = z^p}

\title{Modular congruences, Q-curves, and the diophantine equation $x^4 + y^4 = z^p $}
\author{Luis V. Dieulefait  \thanks{supported by
Post-doctoral Fellowship from the European Research Network
``Galois Theory and Explicit Methods
 in Arithmetic"  at the Institut de Math{\'e}matiques de Jussieu
  and by a MECD postdoctoral grant at the Centre de Recerca Matem{\'a}tica
from Ministerio de Educaci{\'o}n y Cultura}\\
Centre de Recerca Matem{\'a}tica \\  Apartat 50, E-08193 Bellaterra, Spain\\
e-mail: LDieulefait@crm.es}

\begin{document}

\maketitle
\begin{abstract}
We prove two results concerning the generalized Fermat equation $\eq$. In
particular we prove that the First Case is true if $p \neq 7$.
\end{abstract}

\section{Introduction}
In this note we will prove the following results concerning the
generalized Fermat equation $x^4+y^4= z^p$:

\begin{teo}
\label{teo:primo} Let $p$ be a prime such that $p \not\equiv -1 \pmod{8}$
and $p>13$. Then the diophantine equation $\eq$ has no solutions
$x,y,z$ with $(x,y)=1$ and $x y \neq 0$.
\end{teo}
We will call a solution primitive if $(x,y)=1$, and non-trivial if $x y \neq
0$.

\begin{defin}
\label{teo:firstcase} Borrowing the terminology introduced by
Sophie-Germain in connection with Fermat's Last Theorem, we say that a
primitive solution $(x,y,z)$ of $\eq$ is in the First Case if $p \nmid x y$.
\end{defin}

\begin{teo}
\label{teo:secondo} Let $p$ be a prime different from $7$. Then
the diophantine equation $\eq$ has no solutions in the First Case.
\end{teo}

First of all, we have to stress that we will depend heavily on the work
of Ellenberg (see [E]) on the more general equation $x^2 + y^4 = z^p$.
The reader should keep in mind that the following much stronger result
is proved in [E]:
\begin{teo}
\label{teo:Jordan} Let $p$ be a prime, $p \geq 211$. Then the equation
$x^2 + y^4 = z^p$ has no primitive solutions in non-zero integers.
\end{teo}

Our results are not including in this theorem only for the lower bound $p \geq
211$.\\
In [E], a Q-curve $E$ of degree $2$ defined over $\Q(i)$ is attached to a given
non-trivial solution of $x^2 + y^4 = z^p$, and using the modularity of
$E$ (proved in [ES]) it is shown that there is a congruence modulo $p$
between the modular form corresponding to $E$ and a modular form of
weight $2$, trivial nebentypus, and level $32$ or $256$. All cusp forms
in these spaces have complex multiplication (CM). This implies that the
prime $p$ is dihedral for the $\mod \; p$ Galois representation attached to $E$ (this is a
representation of the full Galois group of $\Q$, see [ES], [E] for precise
definitions). To prove theorem \ref{teo:Jordan} in [E], it is proved
(using in particular some results on the Birch and Swinnerton-Dyer conjecture
and  analytic estimates for zeros of special values of L-functions)
that for $p \geq 211$ if the $\mod \; p$ representation
corresponding to a degree $2$ Q-curve over $\Q(i)$ is dihedral then
 the curve must
 have potentially
 good reduction at every prime with residual characteristic greater than $3$.\\
We will follow a different path in order
 to obtain a result holding also for small primes: we
will use the theory of modular congruences (level raising results of Ribet)
 to study the congruences
between the modular form associated to
$E$ and the particular CM modular forms of level $32$ and
$256$. We will restrict to the simpler case where the curve $E$ is
associated to a solution of the equation $\eq$ to obtain a stronger
result (if we start with a solution of $x^2+y^4 = z^p$ an imitation of the arguments
in this paper  only proves the
First Case for primes $p \equiv 1,3 \pmod{8}$, $p>13$, and the non-existence of
 solutions as in theorem \ref{teo:primo}
  only for $p \equiv 1 \pmod{8}$, $p >13$).\\
The third main ingredient in the proof  is the theory of sum of two
squares (results of Fermat) and its relation with the CM cusp form of
level $32$. Let us explain this relation before getting into the proof
of theorems \ref{teo:primo} and \ref{teo:secondo}:
\\ Let $f_1$ be the cusp form in $S_2(32)$. Using the fact that $f_1$
is CM and it verifies $f \cong f \otimes \psi$ where $\psi$ is the
Dirichlet character of conductor $4$, we know that $f_1$ is induced
from a Hecke character of $\Q(i)$. From this, we easily get the
well-known relation for the Hecke eigenvalues $\{ a_q \}$ of $f_1$:\\
$ a_q = (\alpha + \beta i) + (\alpha - \beta i) = 2 \alpha$, with
$\alpha^2 + \beta^2 = q$, if the prime $q$ verifies $q \equiv 1
\pmod{4}$, and $a_q= 0$ if $q \equiv 2,3 \pmod{4}$.\\
Acknowledgement: I want to thank Virginia Balzano.

\section{The results we need}
Let $A,B,C$ be a primitive solution of $A^4 + B^4 = C^p$ with $A B \neq
0$. It is an elementary exercise to see that $(6,C)=1$. Thus, we can
assume that $A$ is even. Following Darmon and Ellenberg, the following
two elliptic curves $E_{A,B}$ and $E_{B,A}$ can be attached to this
triple:
$$ E_{A,B}: \quad \quad y^2 = x^3 + 2(1+i) A x^2 + (-B^2 + i A^2) x $$
$$ E_{B,A}: \quad \quad y^2 = x^3 + 2 (1+i) B x^2 + (A^2 + i B^2) x $$
They are both degree $2$ Q-curves, i.e., each of them is isogenous
 to its Galois conjugate, with a
degree $2$ isogeny. Both have good reduction at primes not dividing $2
C$, and because $3 \nmid C$ this already implies (cf. [ES]) that they
are modular. Let us denote by $E$ any of these two Q-curves, whenever we
do not need to distinguish between them.\\
Modularity should be interpreted in terms of the compatible family of
Galois representations of $G_\Q$ attached to $E$ (see [ES] for
definitions):
$$ \rho_{E,\lambda}: G_\Q \rightarrow \GL_2(\Q(\sqrt{2})_\lambda)$$
for $\lambda$ a prime in $\mathbb{Z}(\sqrt{2})$. Each $\rho_{E,\lambda} $
 is unramified
 outside $2 C \ell$, $\ell$ the rational prime below $\lambda$.
 These representations are modular and by construction they correspond
 to a modular form $f$ with $\Q_f = \Q(\sqrt{2})$ having an extra twist
 given by the character $\psi$ corresponding to $\Q(i)$: $f^{\sigma} \cong
 f \otimes \psi$, $\sigma$ generating $\Gal(\Q(\sqrt{2})/\Q)$.\\
In [E], generalizing results of Mazur to the case of Q-curves, it is
shown that if $\ell > 13$ the residual representation
$\bar{\rho}_{E,\lambda}$ is irreducible.\\
From now on we will assume $p> 13$. The close relation between the
discriminant of $E$ and $C^p$ shows that for $P \mid p$ in
$\mathbb{Z}(\sqrt{2})$ the residual $\mod \; P$ representation $\bar{\rho}_{E,P}$
 has conductor equal to a power of $2$. The exact value of this
 conductor was computed in [E], giving $32$ for the case of $E_{A,B}$ (recall $A$ is
  even) and $256$ for $E_{B,A}$.\\
 The modularity of both Q-curves together with Ribet's level-lowering
 result give:
 $$ \bar{\rho}_{E_{A,B}, P} \cong \bar{\rho}_{f, P} \quad \quad f \in S_2^{new}(32)$$
 $$ \bar{\rho}_{E_{B,A}, P} \cong \bar{\rho}_{f', P} \quad
  \quad f' \in S_2^{new}(256)$$
All newforms of these levels have CM. Thus, this implies that the
(projective) images of $ \bar{\rho}_{E_{A,B}, P} $ and $\bar{\rho}_{E_{B,A}, P}$ fall
both in the normalizer of a Cartan subgroup of $\PGL_2(\mathbb{F}_p)$.
 Allways with the assumption $p >
13$, generalizing results of Momose
 to the case of Q-curves, it is
proved in [E] that the case of a split Cartan subgroup is impossible,
the case of a non-split Cartan subgroup remaining the only
case to be considered.

\section{The Proofs}
Let us start by describing in more detail the newforms of levels $32$
and $256$:\\
$N=32$: The cusp form $f_1 \in S_2(32)$ corresponds to the elliptic
curve $y^2 = x^3-x$ which has CM by $\Q(i)$.\\
$N=256$: $f_2$ (and its Galois conjugate) correspond to the degree $2$ Q-curve
$y^2 = x^3 + 2 (1+i) x^2 + i x$ which has CM by $\Q(\sqrt{-2})$. We
have $\Q_{f_2} = \Q(\sqrt{2})$, $f_2$ has thus both CM and an inner
twist (which is not uniquely defined).\\
$f_3$ and $f_4$: Corresponding to elliptic curves defined over $\Q$
with CM by $\Q(i)$.\\
$f_5$ and $f_6$:  Corresponding to elliptic curves defined over $\Q$
with CM by $\Q(\sqrt{-2})$.\\
The eigenvalues $a_p$, $p \leq 17$, of these cusp forms, are:\\
$$f_1: 0,0,-2,0,0,6,2$$
$$f_2: 0, 2\sqrt{2}, 0, 0 , -2\sqrt{2}, 0 ,6$$
$$ f_3: 0,0,-4,0,0,-4,-2$$
$$f_4: 0,0,4,0,0,4,-2$$
$$f_5: 0,-2,0,0,-6,0,-6$$
$$f_6: 0,2,0,0,6,0,-6$$
The first thing to observe is that $f_5$ and $f_6$ can be eliminated
from the possibilities (this result will not be essential in the sequel but
it may be of independent interest):
\begin{prop}
\label{teo:cincoyseis} Let $E$ be a Q-curve of degree $2$ defined over
$\Q(i)$ and with good reduction  at $3$. Let $p>13$ be a prime and $P \mid p$
such that: $\bar{\rho}_{E,P} \cong \bar{\rho}_{f_t,P}$
for $t \leq 6$, where $f_t$ is one of the cusp forms described above.
Then $t \leq 4$.
\end{prop}
Proof: We know that $\rho_{E,P}$ is defined over $\Q(\sqrt{2})$ and has
an extra twist: $\rho_{E,P}^\sigma \cong \rho_{E,P} \otimes \psi$
where $\psi$ is the $\mod \; 4$ character. This implies that for every
good reduction prime
$q \equiv 3 \pmod{4}$: $a_q = z \sqrt{2}$, $z$ a rational integer;
where $a_q = \trace (\rho_{E_P}(\Frob \; q))$, $P \nmid q$.\\
Using the bound $|a_3| \leq 2 \sqrt{3}$ we have the only possibilities:
$a_3 = 0, \pm \sqrt{2}, \pm 2\sqrt{2} $. The cusp forms $f_5$ and $f_6$
have Hecke eigenvalue $a_3 = \pm 2 $, so the congruence
$\bar{\rho}_{E,P} \cong \bar{\rho}_{f_t , P}$ with $t=5,6$ gives at
$q=3$: $a \equiv \pm 2 \pmod{P}$ for $a \in \{ 0, \pm \sqrt{2} , \pm 2 \sqrt{2}
\}$, but this is impossible for $p>2$.\\

Proof of Theorem \ref{teo:primo}:\\
We have attached two Q-curves to a
 primitive non-trivial solution $A,B,C$ of $A^4+B^4=C^p$, and for the Galois
 representations attached to them we have the congruences:
 $$ \bar{\rho}_{E_{A,B}, P} \cong \bar{\rho}_{f_1, P} \quad \quad (3.1)$$
 $$ \bar{\rho}_{E_{B,A}, P} \cong \bar{\rho}_{f_t, P} \quad \quad (3.2)$$
 for $t=2,3$ or $4$.\\
 Following [E], for $p>13$ we can assume that the projective images lie
 both in the normalizer  of non-split Cartan subgroups of
 $\PGL_2(\mathbb{F}_p)$. Using the fact that $f_1$ has CM  by $\Q(i)$,
 we know that for $E_{A,B}$ this is the case precisely when $p \equiv 3
 \pmod{4}$. This proves the theorem for $p \equiv 1 \pmod{4}$,
 $p>13$.\\
 So let $p$ be a prime $p \equiv 3 \pmod{4}$, $p >13$. Observe that
 this implies in particular that $E_{A,B}$ and $E_{B,A}$ have both good
 reduction at $p$, because from the equation $A^4+B^4=C^p$ and
 $(A,B)=1$ it is a very old result that all primes dividing $C$ are of the
 form $4k+1$. Let $\{a_q \}$ be the set traces of $\rho_{E_{A,B}, P}$
 for $q \nmid 2 C p$. Congruence (3.1) gives
 at $q=3$: $
 a_3 \equiv 0 \pmod{P}$. But we know that $a_3 \in  \{ 0, \pm \sqrt{2}, \pm 2 \sqrt{2}
 \}$, then the only possibility is $a_3 = 0$. The
  value $a_3$ is
 computed as usual by counting the number of points of the reduction of
 the curve $E_{A,B}$ modulo $\check{3} $ for a prime $\check{3}$ dividing $3$.
 It depends only on the values of $A$ and $B^2$ modulo $3$. A direct
 computation shows that $a_3 = 0$ if and only if $A \equiv 0 \pmod{3}$
 and $B \not\equiv 0 \pmod{3}$. Thus, we can assume that $A$ and $B$
 verify these congruences (recall that by assumption $A$ is even, so we
 have in particular $6\mid A$). This determines up to sign the value of the trace
 $a'_{3}$ of the image of $\Frob \; 3$ for the Galois representation $\rho_{E_{B,A},P}$ (it depends
  only on the values of $A$ and $B$ modulo $3$). In fact, putting $A=0$ and
  $B = \pm 1$
  we obtain  $a'_{3} = \pm 2 \sqrt{2}$. \\
  If we suppose that congruence (3.2) holds for $t=3$ or $4$, comparing
  traces at $q=3$ we obtain: $a'_3 = \pm 2 \sqrt{2} \equiv 0 \pmod{P}$,
  and this is a contradiction. So we conclude that if $B,A,C$ is a
  non-trivial primitive solution of $B^4+A^4=C^p$ with $A$ even and
  $p>13$, then congruence (3.2) has to be verified by the newform
  $f_2$ (or its Galois conjugated). Applying again Ellenberg's generalization
  of  results of Momose (cf. [E]) to the case of Q-curves, we can
  assume that the projective images
   of the congruent $\mod \; P$ representations in (3.2)
   must be contained in the normalizer of a non-split Cartan subgroup,
   and using the fact that $f_2$ has CM by $\Q(\sqrt{-2})$ we know that
   this can only happen if $p \equiv 5,7 \pmod{8}$. This proves the
   result also for $p \equiv 3 \pmod{8}$, $p>13$, which concludes the
   proof.\\

   Proof of Theorem \ref{teo:secondo}:\\
   Assume as before that $p>13$ and consider congruence (3.1) again.
   As we already mentioned in the proof of theorem \ref{teo:primo},
    we can restrict to $p \equiv 3
   \pmod{4}$ (the non-split Cartan case).
   The odd primes where $E_{A,B}$ has bad reduction are
   precisely the primes $q \mid C$; the curve
   has semistable reduction at these primes (cf. [E]) and we obviously
   have $q \equiv 1 \pmod{4}$. In particular $E_{A,B}$ has good
   reduction at $p$.\\
   As explained in [E], from the fact that $E_{A,B}$ has multiplicative
   reduction at $q \nmid p$, and the assumption that the corresponding
   projective $\mod \; P$
   Galois representation has image in the normalizer of a non-split
   Cartan subgroup (using the fact that the cusps of $X_0^{ns}(2,p)$ have minimal
    field of definition $\Q(\zeta_p + \zeta_p^{-1})$, with $\zeta_p$ a primitive
    $p$-th root of unity) it follows that the residue field $\mathbb{Z}[i]/\check{q}$
    must contain $ \zeta_p + \zeta_p^{-1} $ for $ \check{q} \mid q $,
    and this implies (because $q$ splits in $\Q(i)$) that $q^2 \equiv 1
    \pmod{p} \quad \; \; (3.3)$. \\
    On the other hand, the level of the modular form $f$ corresponding to $E_{A,B}$
    is $32 \cdot \cond(C)$, where $\cond(C)$ is the product of the primes
    dividing $C$. This implies, by Ribet's level-lowering result, that
    together with congruence (3.1) there are congruences:
    $$ \rho_{f,P} \equiv \rho_{f_1, P} \equiv \rho_{f_{(q)}, P} \quad \pmod{P}$$
    for every prime $q \mid C$, with $f_{(q)}$ a newform of level $32
    q$. Ribet's level raising result gives a constraint for such
    congruence primes (see [G]): using the fact that $f_1$ has level $32$ and
    $f_{(q)}$ has level $32 q$, the above congruences imply that:
    $$ a_q^2 \equiv q (q+1)^2 \pmod{P}$$
    where $a_q$ is the $q$-th coefficient of $f_1$, for every $q \mid
    C$.\\
     Using (3.3),  in the above congruence we substitute $q$ by $\pm 1$
     and we obtain: $a_q^2 \equiv 4,0 \pmod{P} $ (convention:
     we list first the value corresponding
     to $q \equiv 1$).\\
    As we explained in section 1, we have $a_q = 2 \alpha$ with
    $\alpha^2 + \beta^2 = q$. Thus, the above congruence gives:
    $\alpha^2 \equiv 1,0 \pmod{p} \quad (3.4)$.\\
    In particular, if $q \equiv -1 \pmod{p}$ we have: $\alpha \equiv 0
    \pmod{p}$, then $\beta^2 = q - \alpha^2 \equiv -1 \pmod{p}$, but
    this is impossible because $p \equiv 3 \pmod{4}$.\\
    Then, it must hold $q \equiv 1 \pmod{p}$ and from (3.4): $\alpha^2 \equiv 1
    \pmod{p}$. Thus, $\beta^2 = q- \alpha^2 \equiv 0 \pmod{p}$.\\
    We have proved that for every $q \mid C$, if we write it as a sum
    of squares
     $q= \alpha^2 + \beta^2$, these two integers
    have to verify: $\alpha^2 \equiv 1 \pmod{p}$ and $p \mid \beta$ (or viceversa).
    Applying the product formula we see that if we write $C^p = R^2 +
    S^2$ in any possible way, it must always hold $p \mid R S$. \\
    Thus, in the equation $A^4 + B^4 = C^p $ we have $p \mid AB$, and
    this proves the theorem for $p>13$. The First Case for the
    remaining small primes different from $7$ (and more generally for $p \not \equiv
    \pm 1 \pmod{8}$) was already solved  … la Kummer
     (see [P] and [C]).

\section{Bibliography}

[C] Cao, Z. F.,{\it The diophantine equation $c x^4 + d y^4 = z^p$}, C. R.
 Math. Rep. Acad. Sci. Canada, 14 (1992), 231-234\\

[E] Ellenberg, J.,{\it Galois representations attached to Q-curves
and the generalized Fermat equation $A^4 + B^2 = C^p$}, preprint\\

[ES] Ellenberg, J., Skinner, C., {\it On the modularity of Q-curves},
Duke Math. J., 109 (2001), 97-122\\

[G] Ghate, E., {\it An introduction to congruences between modular
forms}, in
``Current Trends in Number Theory", proceedings of the International
 Conference on Number Theory,
Harish-Chandra Research Institute, Allahabad, November 2000;
S. D. Adhikari, S. A. Katre, B. Ramakrishnan (Eds.), Hindustan Book
Agency (2002) \\

[P] Powell, B.,{\it Sur l'{\'e}quation diophantienne $x^4 \pm y^4= z^p$}, Bull.
 Sc. Math., 107 (1983), 219-223\\

\end{document}